\documentclass[a4paper, 10pt]{amsproc}

\usepackage{defs}
\usepackage{amsmath}
\usepackage{fancyvrb}
\usepackage{graphicx}
\usepackage{tabularray}
\usepackage{adjustbox}
\usepackage{booktabs}
\usepackage{listings}
\usepackage{xcolor}
\usepackage{orcidlink}
\usepackage{tikz}
\usetikzlibrary{shapes.geometric, arrows}
\usepackage[linesnumbered,ruled,noline,slide,vlined,commentsnumbered]{algorithm2e}
\SetKwRepeat{Do}{do}{while}%
\usepackage[noend]{algpseudocode}

\usepackage{orcidlink}
\hypersetup{final=true}
\lstdefinestyle{Matlab-editor}{
    language=Matlab,                 
    basicstyle=\ttfamily\small,      
    keywordstyle=\color{blue},       
    commentstyle=\color{green!60!black}, 
    stringstyle=\color{orange},      
    numbers=left,                    
    numberstyle=\tiny\color{gray},   
    stepnumber=1,                    
    frame=single,                    
    rulecolor=\color{gray},          
    backgroundcolor=\color{white},   
    showspaces=false,                
    showstringspaces=false,          
    breaklines=true,                 
}

\title[Walker slice sampling based global optimization algorithm]{On a probabilistic global optimizer derived from the Walker slice sampling}

\author[A. Gupta]{Aditya Gupta\orcidlink{0009-0000-7121-1284}}
\author[S. Das]{Souvik Das\orcidlink{0000-0001-6918-6219}}
\author[D. Chatterjee]{Debasish Chatterjee\,\orcidlink{0000-0002-1718-653X}}

\begin{document}

\maketitle

\begin{abstract}
	This article presents a zeroth order probabilistic global optimization algorithm --- \(\algoname\) --- for (not necessarily convex) functions over a compact domain. A discretization procedure is deployed on the compact domain, starting with a small step-size \(h > 0\) and subsequently adaptively refining it in the course of a simulated annealing routine utilizing the Walker slice and the Gibbs sampler, in order to identify a set of global optimizers up to good precision. \(\algoname\) is parallelizable, which helps with scalability as the dimension of decision variables increases. Several numerical experiments are included here to demonstrate the effectiveness and accuracy of \(\algoname\) in high-dimensional benchmark optimization problems.
\end{abstract}

\section{Motivation and significance}
Global optimization is a challenging topic, and a bewildering array of techniques have been invented to solve global optimization problems that arise in the natural and engineering sciences over the past several decades. Each technique is either designed to tackle specific classes of problems, or tends to work best for special classes of global optimization problems; our contributions in this article follow this trend closely, being motivated by peculiarities that arise in global optimization problems in the context of a recent algorithm (reported in \cite{ref:DasAraCheCha-22}) for optimal solutions to convex semi-infinite optimization problems. The ubiquity of this class of robust convex optimization problems, together with the indicated new optimal algorithmic technique whose performance depends on the fidelity of a global optimization routine, are the chief reasons for our engagement in this direction.

Convex semi-infinite programs (CSIPs) originate in a diverse array of robust convex optimization problems that appear naturally in a range of applications; we refer the reader to the review articles \cite{RH-OKK:93,MAL-GS:2007,MAG-MAL:18} for more information. CSIPs feature a convex objective function to be minimized over the intersection of a compact family of convex constraint sets, and since this family of constraints may be uncountably large, it poses a stiff challenge for numerical techniques designed to solve CSIPs. CSIPs are NP-hard \cite[Chapter \(16\)]{ref:AA-IT-21}, and algorithmic techniques for solving CSIPs until \(2022\) were of the infinitary type --- either relying on infinite memory or infinite computational power for exact solutions in general. Moreover, quantitative error estimates after finite truncation were, in general, unavailable.

Relying on a structural result in \cite{ref:Bor-81}, a new computationally tractable algorithm to extract optimal solutions to CSIPs was established in \cite{ref:DasAraCheCha-22}. It relies on an iterated max-min optimization problem, of which the outer maximization is a finite-dimensional (zeroth order) global optimization that scales \emph{linearly} with the dimension of the decision space of the original CSIP. Theoretically, this algorithm, which we shall refer to as the MSA algorithm hereinafter,\footnote{The term MSA derives from the names of three students, Mishal Assif P. K., Souvik Das, and Ashwin Aravind, who contributed to its development.} provides exact solutions to CSIPs, but in practice, one must contend with the challenges posed by the accurate computation of solutions to the global optimization stage therein for which the value of the objective can be evaluated at will, but admits no formula. The current article takes the first steps toward developing a dedicated probabilistic zeroth order global optimization solver for problems that naturally arise in the process of solving CSIPs via the MSA algorithm. 


\subsection*{Our Contributions}
\label{sec: 2}
Against the preceding backdrop, the key advancements are:
\begin{enumerate}[label=(\Alph*), leftmargin=*, widest=b, align=left]
    \item \textbf{A zeroth order global optimization solver:} Our zeroth order solver \(\algoname\) is of independent interest in the context of global optimization, and is naturally relevant in the broader areas of machine learning, data science, and control theory due to its seamless applicability to optimization problems featuring \emph{non-convex} objective functions over compact domains. \(\algoname\) does not impose strict regularity assumptions on the objective function such as differentiability or smoothness; \(\algoname\) is fine-tuned to high-dimensional problems, and are easily parallelizable, all of which are, we believe, crucial to several data science applications. 

    \item \textbf{Key features:} \(\algoname\) has been designed for solving CSIPs via the MSA algorithm, and in this context, the key features of the simulated annealing routine in \(\algoname\) are:
    \begin{itemize}[label = \(\circ\), leftmargin=*]
		\item The compact domain is discretized progressively finely after starting from a sufficiently small discretization step \(h > 0\). This initial phase aids in quick exploration, and if certain conditions —-- elaborated in Algorithm \ref{alg:adapref} --- are satisfied, then the discretization step size \(h\) is refined (reduced) \emph{adaptively}; see Algorithm \ref{alg:adapref} in Appendix \ref{appen:details} for more details, resulting in better accuracy without adversely affecting its speed.

		\item  The optimization in \(\algoname\) itself involves simulated annealing on the discretized domain via Markov chain Monte Carlo with the \emph{Walker slice} sampling algorithm \cite{ref:SGW-14} along each dimension, along with the adaptive-refinement scheme mentioned in the preceding point. The Walker slice sampling facilitates quicker sampling relative to the more common Metropolis-Hasting scheme, and Gibbs sampling was incorporated to sample across multiple dimensions simultaneously in parallel; the readers are referred to \cite{Kirkpatrick1984,Cerny1985} for more details on the simulated annealing routine and \cite{ref:Hagg-02} for a brief exposition on the Gibbs sampling strategy.
    \end{itemize}
\end{enumerate}
Extensive numerical experiments were conducted with \(\algoname\). It was tested on two benchmark test functions known as the Ackley and Levy function, in high dimensions; compared to the recent algorithm ProGO reported in \cite{ref:Zhang23}, we observed a substantial improvement in the accuracy at the cost of longer simulation horizons of \(\algoname\) --- see \S \ref{sec:4} for more details. We tested \(\algoname\) on two benchmark optimization problems lifted from the literature and the results are reported in Tables \ref{tab:minlp} in \S \ref{sec:4}. Moreover, to demonstrate the effectiveness of \(\algoname\) in solving CSIPs via the MSA algorithm, \(\algoname\) was applied to a mid-scale semi-definite optimization problem lifted from \cite{ref:BB-99}, and a particularly difficult instance of the Chebyshev center problem reported in \cite{ref:PP-DC-23}. The latter is known to be an NP-hard problem, and to date, approximate algorithms focusing only on relatively simpler settings have been reported in the literature --- see \cite[Chapter \(15\)]{ref:AA-IT-21} for more details. \(\algoname\) computes the Chebyshev radius and the Chebyshev center up to good accuracy compared to the one reported in \cite{ref:PP-DC-23} (where the authors employed a deterministic simplicial technique to get the optimal value and optimizers).

\subsection*{Notations used in the article}
For us \(\N \Let \aset{1,2,\ldots}\) is the set of positive integers, \(\Nz\Let \N \cup \aset[]{0}\), and \(\Z\) is the set of integers. The \(d\times d\)-dimensional identity and zero matrices are \(\identity{d} \in \Rbb^{d \times d}\) and \(\zv{d} \in \Rbb^{d \times d}\), respectively. All random variables are assumed to be defined on a fixed probability space \((\Omega,\tsigalg, \PP)\). The notation \(x \sim \textsf{Unif}([a,b])\) means that \(x\) is uniformly distributed over the interval \([a,b]\). Moreover, we denote by \(|V|\) the cardinality of the set \(V\).

\subsection*{Organization:}
The article unfolds as follows: \secref{sec:3} introduces our algorithm \(\algoname\) (Algorithm \ref{alg:swiftnav}) with the relevant background and the associated technical assumptions crucial for \(\algoname\). In \secref{sec:4} we present the experimental results conducted to demonstrate the fidelity of \(\algoname\), which is followed by Conclusion in \S\ref{sec:conclu_dissc}.

\section{Problem description}\label{Sec:SwiftNav_PS}
This article presents an algorithm \(\algoname\) to solve global optimization problems\footnote{A detailed description along with the requisite codes can be found in this GitHub repository: \href{https://github.com/souvik84/SwiftNav_V1}{\texttt{https://github.com/souvik84/SwiftNav\_V1}.}}
\begin{equation}\label{eq:global opt prob}
    \begin{aligned}
        & \inf_{x \in \admst} \; \obj(x)
    \end{aligned}
\end{equation}
with the following data:
\begin{enumerate}[label=(\ref{eq:global opt prob}-\Alph*), leftmargin=*, widest=b, align=left]
    \item \label{it:objective function} \(f: \admst \subset \Rbb^n \lra \Rbb \) is the (not necessarily convex) objective function, assumed to be measurable;

    \item \label{it:domain} \(\admst \neq \emptyset\) is a compact domain subset of \(\Rbb^n\), but not necessarily convex. 
\end{enumerate}
We assume that the optimization problem \eqref{eq:global opt prob} admits a global solution; i.e., the set of optimizers is nonempty. In what follows in \(\algoname\), we uniformly discretize the compact set \(\admst\) with a step size \(h\) (which may be refined adaptively depending upon the desired accuracy), and the discretized version of \(\admst\) is denoted by \[\stsp \Let \prod_{i=1}^n \stsp_i \subset h \bigl(\underbrace{\Zbb \times \cdots \times \Zbb}_{n\text{-times}} \bigr).\]
The rest of this article focuses on different aspects and functionalities of \(\algoname\). 

\section{\(\algoname\): Our algorithm}
\label{sec:3}
We adopt the established notations in \S\ref{Sec:SwiftNav_PS}, and denote by \(\stsp_i\) the state-space of the Markov chain along the \(i^{\text{th}}\) dimension.
\subsection{Architecture}
\label{subsec:swiftnav details}
\begin{figure}[!htpb]
    \centering
     \includegraphics[scale=0.4]{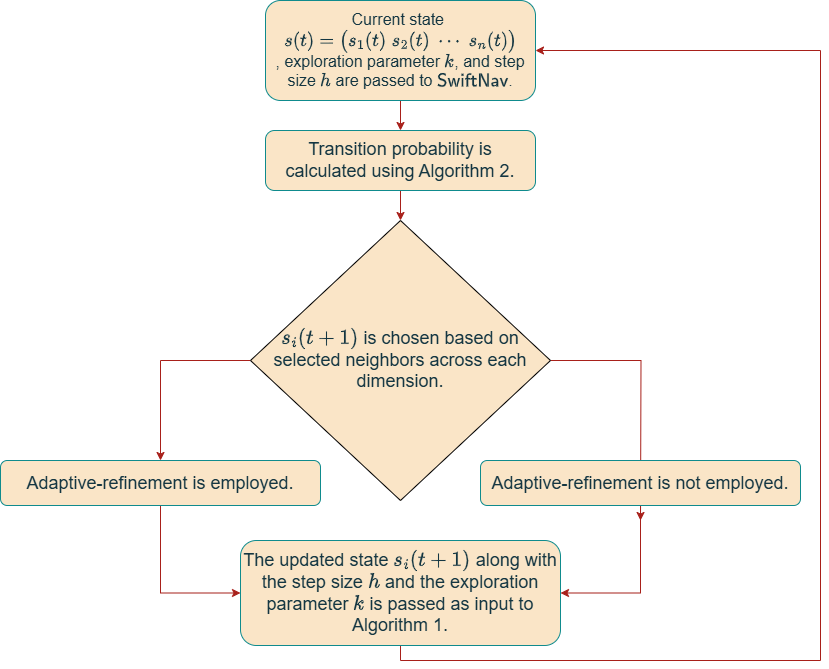} 
     \caption{Schematic of \(\algoname\).}
    \label{fig:flow}
\end{figure}

With discretization step size \(h>0\) picked as above, we perform simulated annealing with Monte-Carlo Markov chain (MCMC) between successive cooling times via the Walker slice sampler in conjunction with the Gibbs technique, one dimension at a time. Independent sampling across dimensions allows our algorithm to be parallelizable (depending on availability, over \(n\) processors),  enabling the simultaneous use of multiple CPUs, with one dedicated to each dimension.\footnote{Parallel processing in Algorithm \ref{alg:swiftnav} is employed using MATLAB's \textit{parfor} loop iterator.} The discretization step is refined iteratively after the preceding annealing routine reports reasonable convergence; the refinement stops if there is no improvement in the objective function after further refinement.
A schematic of Algorithm \ref{alg:swiftnav} from a bird's-eye view is provided in  Figure \ref{fig:flow}. The notations related to Algorithm \ref{alg:swiftnav} are relegated to \S\ref{subsec:Notation Swiftnav}.
\subsection{Key functionalities}
\label{subsec:Soft_func}
For the \(i\)$^{\text{th}}$ dimension and for each time \(t\), let us fix a state \(s_i(t) \in  \stsp_i\) of the underlying Markov chain. Define by \(\neigh_i(t) \subset h\Zbb\) the set neighboring points of \(s_i(t)\) and let \(\neigh(t) \Let \bigcup_{i=1}^n \neigh_i(t)\). Let \(s(t+1) = \bigl(s_1(t),s_2(t),\ldots, s_n(t) \bigr)\)  and we employ the notation \(s_{-i}(t) = \bigl(s_1(t),s_2(t),\ldots,s_{i-1}(t),s_{i+1}(t),\ldots, s_n(t) \bigr) \).
At the time step \(t + 1\):
\begin{enumerate}[leftmargin=*,label=(\Alph*)]
 \item \label{it:pointA} \(\ol{r}_i(t) \in \neigh_i(t)\) is sampled according to \(p_A\cprobof[\big]{\ol{r}_i(t) \given s_{-i}(t)}\) along the \(i\)$^{\text{th}}$ dimension, utilizing a Gibbs sampling scheme \cite[Chapter \(7\), pp.no. 49]{ref:Hagg-02} while freezing transitions along all other dimensions;

\item \label{it:pointB}The procedure in \ref{it:pointA} is repeated for each dimension to select a new sample point \(\ol{r}(t) \Let \bigl(\ol{r}_1(t), \ol{r}_2(t), \cdots, \ol{r}_n(t)\bigr) \in \neigh(t).\)
The new vector \(s(t+1) = \bigl(s_1(t+1), s_2(t+1),\cdots, s_n(t+1)
\bigr)\) at time \(t+1\) is then given by \(s(t+1) = \ol{r}(t).\)

\item The steps \ref{it:pointA} and \ref{it:pointB} above are repeated until the Markov chain converges to its stationary distribution.
\end{enumerate}
Let us now expand on the key functionalities of \(\algoname\) in more detail. 
\begin{itemize}[label = \(\circ\),leftmargin=*]
    \item \textbf{Initialization:} \(\algoname\) is initialized at time \(t=0\) with an initial point \(s(0) \sim \textsf{Unif}\bigl( \neigh(0)\bigr)\), where
    \( \neigh(0) \) can either be chosen as the entire domain or a subset of it.

    \item \textbf{Exploration:} Once initialized, the exploration step consists of employing the Walker slice sampling algorithm (see refer to Appendix \ref{subsec:ws} for a brief exposition) to search the domain in the simulated annealing algorithm for a global solution corresponding to the optimization problem \eqref{eq:global opt prob}.  Let us denote the neighbors of \(s_i(t)\) by \(\bigl(\dummyx^m \bigr)_{m=1}^{\abs{V_i(t)}} \subset V_i(t)\), and assume that \(\bigl(\dummyx^m\bigr)_{m=1}^{\abs{V_i(t)}} \) are arranged in ascending order, i.e.,  for a fixed \(t\in \Nz\),
    \(\dummyx^1 \leq \dummyx^2 \leq \cdots \leq \dummyx^{\abs{V_i(t)}}. \)
    At every \(t \in \Nz\), we compute \(p_A \cprobof[\big]{\cdot \given s_{-i}(t)}\)(via Algorithm \ref{alg:probcal} in Appendix \ref{appen:subsec:probcal}) for the neighborhoods \(\bigl(\dummyx^m\bigr)_{m=1}^{\abs{V_i(t)}} \subset V_i(t)\) and
    generate a random number \(\zeta \sim \mathsf{Unif}([0,1])\) such that \(0 < \zeta < 1\) and 
    \begin{equation*}
    \ol{m} \Let  \min \aset[\bigg]{ m \geq 1 \suchthat \sum_{k=1}^{m}p_A\cprobof[\big]{\dummyx^{k} \given s_{-i}(t)} \geq \zeta }.
\end{equation*} 
 Then assign \(\dummyx^{\ol{m}}\) to \(s_i(t+1)\) to be the new state along the \(i^{th}\) dimension, allowing us to choose \(s_i(t+1)\) at time \(t+1\) from  \(V_i(t)\) with a probability of \(p_A\cprobof[\big]{s_i(t+1) \given s_{-i}(t)}\). It is important to note that this procedure is repeated across each dimension separately, with the associated computations being carried out keeping the states for all other dimensions unchanged. The end result is an updated vector \( s(t+1) = \bigl(s_1{(t+1)} \; s_2{(t+1)} \cdots s_n{(t+1)}\bigr)\), and the objective function is evaluated on \(s(t+1)\).

\item \textbf{Adaptive-refinement of the discretization step size \(h\):} The adaptive-refinement scheme is employed during the \emph{exploration} step mentioned above by monitoring the best value of the objective function evaluated at \(r_i(t) \in V_i(t)\). If the algorithm shows no significant improvements after a pre-specified number of iterations, then it refines the discretized space by reducing the step size \(h\) to enhance accuracy. If the adaptive-refinement scheme yields better results, the algorithm further adjusts the step size to improve accuracy. This process is repeated until no further improvement is observed, and the algorithm reverts to the original step size. Readers are referred to Algorithm \ref{alg:adapref} and its associated notations in Appendix \ref{appen:subsec_adapref}, for more details.
\end{itemize}
\subsection{Notation used in Algorithm \ref{alg:swiftnav}}\label{subsec:Notation Swiftnav}
Algorithm \ref{alg:swiftnav} computes the new state \(s(t+1)\) at time \(t+1\) given the current state \(s(t)\). Let us now define the notations used in Algorithm \ref{alg:swiftnav}: 
\begin{enumerate} [leftmargin=*,label=(\Alph*)]
    
    
    \item The variables \texttt{stepSize} and \texttt{stepMatrix} together assist in constructing the matrix `\texttt{neighbours}', whose rows consist of the neighbors of \(s(t)\). 
    The dimension of the matrix `\texttt{neighbours}' is equal to \((2k-1 \times \texttt{numVars})\). Moreover, \(\texttt{neighbours}[\cdot][i]\) denotes the neighbors for the \(i\)$^{\text{th}}$ dimension. Note that `\texttt{neighbors}' is constructed with columns arranged in ascending order to facilitate fast computations.

    

    
    \item  At time \(t\), \texttt{pi\_table} contains the evaluations \(\pi(r_i(t))\) at each neighboring point \(r_i(t)\) of \(s_i(t) \in \neigh_i(t)\) for the \(i^{\text{th}}\) dimension. This step allows Algorithm \ref{alg:probcal} (in Appendix \ref{appen:details}) to avoid repeated evaluations of \(\pi\), thereby, significantly speeding up the process.     
    
    \item \texttt{objectiveFunc()} is used to evaluate the value of the objective function \(\obj\) at \texttt{newState}. This is also passed as an argument to Algorithm \ref{alg:adapref} during the mesh refinement scheme; see Appendix \ref{appen:details} for more details.
    
\end{enumerate}
Algorithm \ref{alg:swiftnav} is compiled as a custom annealing function used with \textsf{MATLAB's} \textsf{simulannealbnd} function and all other settings were kept default.
\begin{algorithm}[htpb]
\SetAlgoLined
\footnotesize
\SetKwInOut{ini}{Initialize}
\SetKwInOut{inp}{Input}
\KwData{exploration parameter \(k\)}
\ini{\(\text{\texttt{currentState}} = s(0) \sim \mathsf{Unif}\bigl([\neigh_1(0) \times \cdots \times \neigh_n(0)]\bigr)\), \(\text{\texttt{numVars}} = n\),  \(\text{\texttt{stepSize}} = h \cdot \text{ones}(\text{\texttt{numVars}},1\)), \(\text{\texttt{stepMatrix}} = (1:k-1)^{\top} \cdot \text{\texttt{stepSize}}^{\top} \)}
\KwResult{\(s(t+1) = \texttt{newState}\)}
\For{\(\texttt{i} = 1:k-1\)}{
    \(\texttt{neighbours(i, :)} \gets \texttt{currentState}^{\top} - \texttt{stepMatrix(i, :)}\)\;
}
\(\texttt{neighbours}(k, :) \gets \texttt{currentState}^{\top}\)\;
\For{\(\texttt{i} = k+1:2k-1\)}{
    \(\texttt{neighbours}(\ttt{i}, :) \gets \texttt{currentState}^{\top} + \texttt{stepMatrix}(\ttt{i}-k, :)\)\;
}
\For{\(\texttt{dim} = 1: \texttt{numVars}\)}{
    \label{line:outer_init}
    \tcp{The outer loop computations in Step \ref{line:outer_init} can be implemented parallelly} 
    \(\texttt{sum} \gets 0\) \;
    \(\texttt{pi\_table} \gets -1\cdot\texttt{ones(2k-1, 1)}\)\; \label{line:pi_init}
    \For{\(\texttt{m} = 1:2k-1 \)}{
        \(\texttt{s} = \texttt{neighbours}(k,\texttt{dim})\) \;
        \(\texttt{r} = \texttt{neighbours}(\texttt{m},\,\texttt{dim}) \) \tcp{\(\ttt{s}\) and \(\ttt{r}\) denotes \(s_i(t)\) and \(r_i(t)\)} \label{line:r_init}
        Compute \(p_A \cprobof[]{\ttt{r} \given \ttt{s}}\) using Algorithm \ref{alg:probcal} for \(\texttt{r}\) in Step \ref{line:r_init} \;
	    \(\texttt{sum} \gets \texttt{sum} + p_A \cprobof[]{\ttt{r} \given \ttt{s}}\) \;
        \(\texttt{nos} \sim \mathsf{Unif}\bigl([0,1]\bigr)\) \;
        \If{\(\texttt{sum} > \texttt{nos}\)}{
            \(\texttt{newState}[m] \gets \texttt{r}\) \tcp{`\texttt{newState}' denotes the updated vector \(s(t+1)\). Choose the neighbour \(\texttt{r}\) as the updated value for the \(\ttt{dim}^{\text{th}}\) dimension}
            \textbf{break}\;
        }
	Update \(\texttt{pi\_table}\) via Algorithm \ref{alg:probcal};
    }
}
Calculate \(\texttt{objectiveFunc(newState)}\)\;\label{line:obj_init}
\(\texttt{currentValue} \gets \texttt{objectiveFunc(newState)}\)\;
Update step size using Algorithm \ref{alg:adapref} \tcp{Refer to \ref{appen:details}}

\caption{\(\algoname\): Exploration of new states}
\label{alg:swiftnav}
\end{algorithm}
%
\section{Illustrative experiments}
\label{sec:4}
\(\algoname\) is a global optimizer and has been evaluated on a range of benchmark problems, both constrained and unconstrained. We assess its performance on high-dimensional Ackley and Levy functions and several constrained optimization problems. Table \ref{tab:compa} and \ref{tab:compl} depict the comparison between \(\algoname\) and the standard Metropolis-Hastings-based simulated annealing algorithm in the context of the Ackley function and the Levy function in $10^3$ dimensions, respectively. 
%
\begin{table}[!htpb]
\centering
\begin{adjustbox}{max width=\textwidth}
\begin{tblr}{lcccccc}
\hline[2pt]
	\SetRow{azure9}
Algorithm & Time taken (seconds) & Log regret & Values attained & Optimal value & Domain \\
\hline[1pt]
\(\algoname\) & 110 & \(\mathbf{0.22997}\) & 1.25856 & 0 & $[-10, 10]^{10^3}$ \\
Simulated Annealing & 110 & \(\mathbf{2.7114}\) & 15.0516 & 0 & $[-10, 10]^{10^3}$ \\
\(\algoname\) & 538 & \(\mathbf{-0.490062}\) & 0.61258 & 0 & $[-10, 10]^{10^3}$ \\
Simulated Annealing & 538 & \(\mathbf{2.7114}\) & 15.0516 & 0 & $[-10, 10]^{10^3}$ \\
\hline[2pt]
\end{tblr}
\end{adjustbox}
\caption{A comparison of \(\algoname\) versus standard Metropolis-Hastings based simulated annealing for the Ackley function in dimension $10^3$.}
\label{tab:compa}
\end{table}
\begin{table}[!htpb]
\centering
\begin{adjustbox}{max width=\textwidth}
\begin{tblr}{lcccccc}
\hline[2pt]
	\SetRow{azure9}
Algorithm & Time taken (seconds) & Log regret & Values attained & Optimal value & Domain \\
\hline[1pt]
\(\algoname\) & 130 & \(\mathbf{2.66177}\) & 14.32161 & 0 & $[-10, 10]^{10^3}$ \\
Simulated Annealing & 130 & \(\mathbf{9.0817}\) & 8793.17 & 0 & $[-10, 10]^{10^3}$ \\
\(\algoname\) & 750 & \(\mathbf{0.65386}\) & 1.92294 & 0 & $[-10, 10]^{10^3}$ \\
Simulated Annealing & 750 & \(\mathbf{8.7038}\) & 6026.2 & 0 & $[-10, 10]^{10^3}$ \\
\hline[2pt]
\end{tblr}
\end{adjustbox}
\caption{A comparison of \(\algoname\) versus standard Metropolis-Hastings based simulated annealing for the Levy function in dimension $10^3$.}
\label{tab:compl}
\end{table}
%
%
%
One of the evaluation metrics used in our experiments is the log regret score: Given $\tilde{x}$ is an estimated optimal solution in a \(d\)-dimensional space, the log regret score is defined by \cite{ref:Zhang23}, for \(f^{\ast} \Let \inf_{x \in \mathbb{X}} f(x),\)
\begin{equation}
r_f = \log \left( f(\tilde{x}) - f^* \right).
\end{equation}
%
\(\algoname\) has been applied to solve convex semi-infinite programs (CSIPs) using the MSA algorithm. Specifically, we integrated the algorithm with the MSA framework to test it on two problems: the Chebyshev center problem for non-convex sets and a semi-definite program (SDP) lifted from \cite{ref:BB-99}. The Chebyshev center problem is challenging \cite[Chapter \(15\)]{ref:AA-IT-21} for non-convex sets owing to the possibly complex geometrical nature of the sets. 
Finding an optimal Chebyshev radius and center for non-convex sets needs the adaptive refinement scheme of \(\algoname\) to capture the behavior of the function \(f\) on pinched regions with significant accuracy relative to, e.g., the native simulated annealing solver reported in \cite{ref:PP-DC-23}.

Algorithm \ref{alg:adapref} in Appendix \ref{appen:details} introduces two parameters, \(p\) and \(g\), which control the number of iterations before which step size is altered. Another parameter, \(\meshfac\), controls the step size \(h\) in the adaptive-refinement scheme. Their values were set at 30, 50, and 2, respectively. The specific values for step size (\(h\)) and exploration parameter (\(k\)) are provided for each case. 

We report that all numerical experiments were conducted on Ubuntu \(14.04.6\) LTS-based server equipped with 64-bit architecture \((\text{x}86\_64)\), \(128\) GB RAM and an Intel Xeon \(\text{E}5\mbox{-}2699 \, \text{v}3\)  processor operating at \(2.30\) Ghz with \(72\) cores.

\subsection*{Code snippets}
The code snippet given below briefly overviews the functionality of \(\algoname\). For illustrations, we take the particular case of the Levy function (the details are provided in \S\ref{subsec:levy function}) as an example.
\begin{enumerate} [leftmargin=*, label=(\Alph*)]
\item \label{it:scrip} The MATLAB script defining the Levy function is added.
\begin{lstlisting}[style=Matlab-editor]
function y = levy(x,n)
   for i = 1:n; z(i) = 1+(x(i)-1)/4; end
   s = sin(pi*z(1))^2;
   for i = 1:n-1
       s = s+(z(i)-1)^2*(1+10*(sin(pi*z(i)+1))^2);
   end 
   y = s+(z(n)-1)^2*(1+(sin(2*pi*z(n)))^2);
end
\end{lstlisting}
\item \label{it:inpp} The algorithm accepts the parameters given in the following table 
via a text file.
\begin{table}[!htpb]
   \centering
   \begin{adjustbox}{max width=\textwidth}
   \begin{tblr}{lccccccc}
   \hline[2pt]
   	\SetRow{azure9}
   Problem Name & Dimensions & \(k\) & \(h\) & Domain & \(p\) & \(f\) & \(q\) \\
   \hline[1pt]
   Levy & \(n\) & 40 & 0.2 & [-10, 10] & 50 & 2 & 30 \\
   \hline[2pt]
   \end{tblr}
   \end{adjustbox}
   \end{table}

\item \label{it:scrip2} \(\algoname\) is implemented as a \emph{custom annealing function}, where a user is provided with a set of options imported from MATLAB's \textsf{simulannealbnd} function, as demonstrated in \ref{it:scrip}. The input file provided in \ref{it:inpp} sets all the parameters for \(\algoname\), mentioned in \ref{it:scrip2}. 
\begin{lstlisting}[style=Matlab-editor]
options = optimoptions(simulannealbnd, AnnealingFcn, {@(x,y) swiftnav(x,y,k,h,p,f,q,false)}, InitialTemperature, initialTemp);
[opt, fval, ~, out] = simulannealbnd(@(x) Ackley(x,n), initialGuess, lb, ub, options);
\end{lstlisting}

\item The algorithm can be implemented on a server or a local computer and the results are observed in real-time via MATLAB's GUI.
\end{enumerate}

\subsection{Ackley function}\label{subsec:ackley}
\begin{figure}[!htpb]
    \centering
     \includegraphics[width=0.6\linewidth]{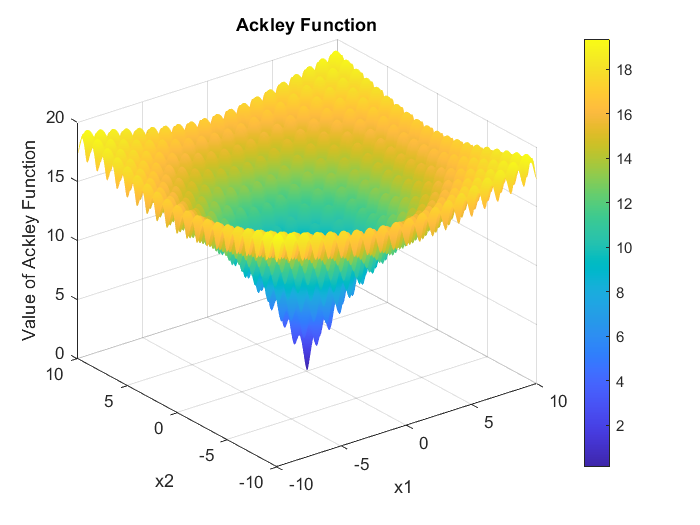} 
     \caption{Plot of the Ackley function on $[-10, 10]^2$.}
    \label{fig:ackplot}
\end{figure}

\begin{figure}[!htpb]
    \centering
     \includegraphics[width=\linewidth]{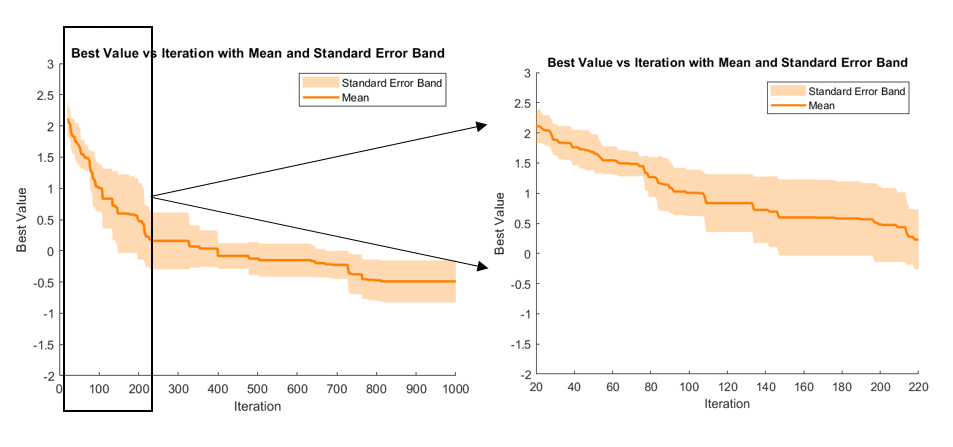} 
     \caption{Plots showing the average log-scaled regret values against iterations for the Ackley function in dimension $10^3$. The curve shows the mean value with standard error represented as a band measured across 10 independent simulations. The complete numerical experiment run for 1000 simulations is on the left and a zoomed-in version showing the first 200 iterations out of 1000 is displayed on the right. Here too, the burn-in time is 20 iterations, which is excluded from the figures.  }
    \label{fig:ackband}
\end{figure}

\begin{figure}[htpb]
    \centering
    \begin{minipage}{0.48\textwidth} 
        \centering
        \includegraphics[width=\linewidth]{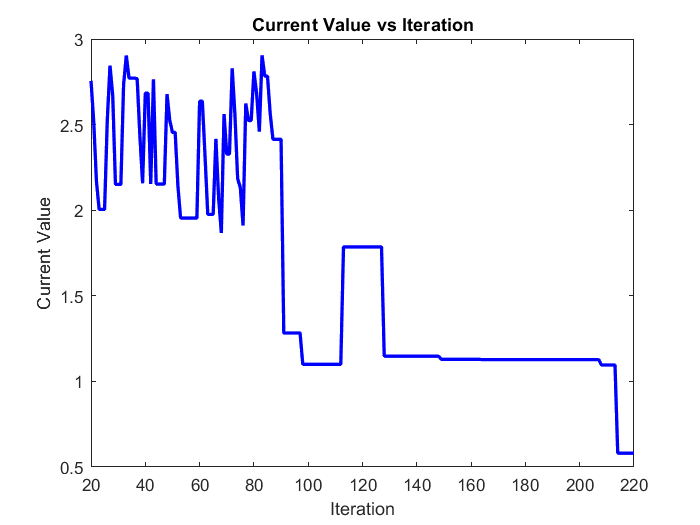} 
        \subcaption{Plot depicting the log regret of the values obtained against iterations.}
        \label{fig:image1}
    \end{minipage}
    \hspace{0.02\textwidth} 
    \begin{minipage}{0.48\textwidth} 
        \centering
        \includegraphics[width=\linewidth]{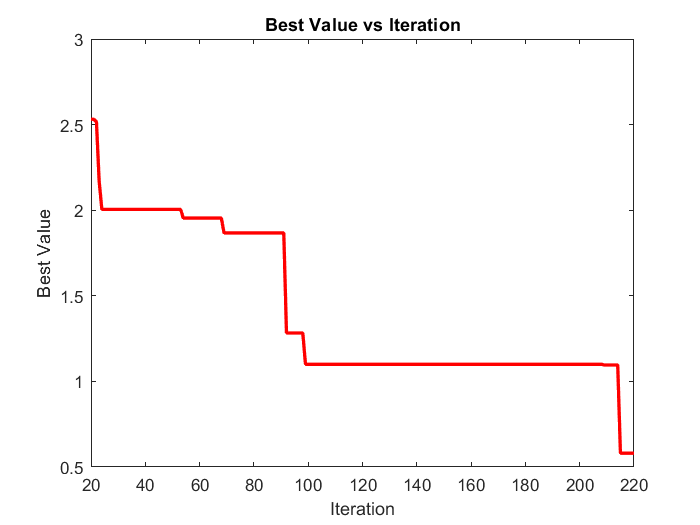} 
        \subcaption{Plot depicting the log regret of the best values obtained thus far against iterations.}
        \label{fig:image2}
    \end{minipage}
    
    \vspace{0.5cm}
    
    \begin{minipage}{0.48\textwidth} 
        \centering
        \includegraphics[width=\linewidth]{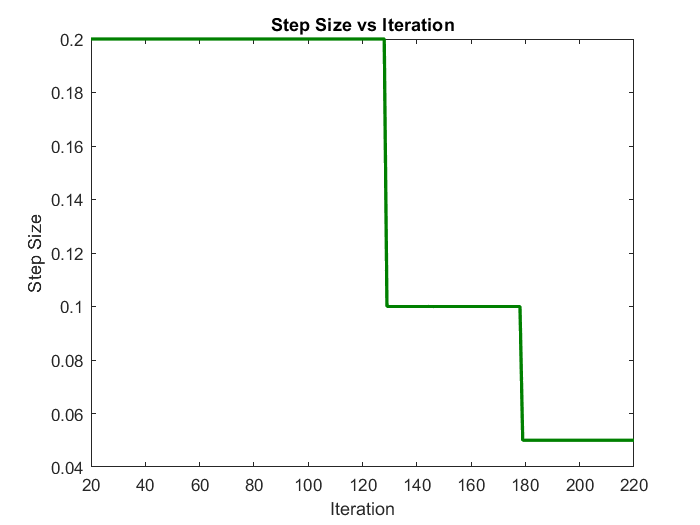} 
        \subcaption{Plot depicting the varying step size trend against iterations.}
        \label{fig:image3}
    \end{minipage}
    \hspace{0.02\textwidth} 
    \begin{minipage}{0.48\textwidth} 
        \centering
        \includegraphics[width=\linewidth]{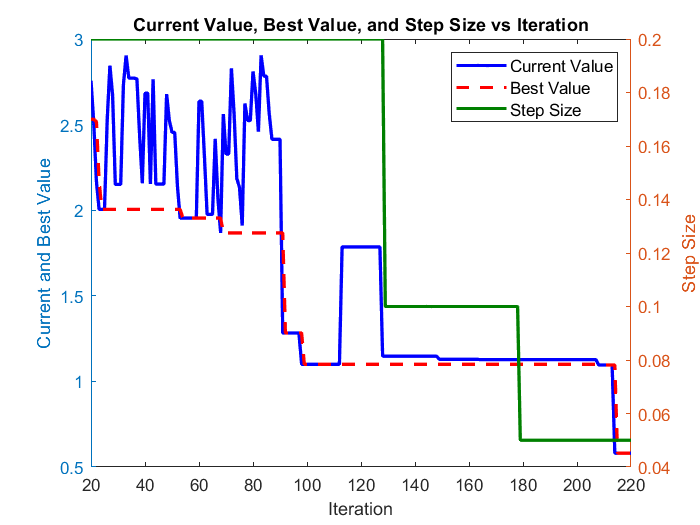} 
        \subcaption{A combined plot showcasing the change in values with the varying step size.}
        \label{fig:image4}
    \end{minipage}
    
    \caption{For the Ackley function the simulation was run with a burn-in time of 20 iterations. Here, we plot the results obtained for the next 200 iterations. Graph (d) displays the changes in values of log regret and step size. A finer mesh can be seen to aid the progression towards global minima.}
    \label{fig:Ack-fig}
\end{figure}

The Ackley function is a widely recognized test function for evaluating the performance of global optimization algorithms. It is known for its complex and multi-modal landscape, which includes numerous local minima. Figure \ref{fig:ackplot} displays a 2-dimensional plot for the Ackley function on $[-10, 10]^2$. The function is defined for an \( n \)-dimensional input vector \( x = (x_1, x_2, \ldots, x_n) \) by the formula:
\begin{equation}
R^n \ni x  \rightarrow f(x) = -a \exp\left( -b \sqrt{ \frac{1}{n} \sum_{i=1}^{n} x_i^2 } \right) - \exp\left( \frac{1}{n} \sum_{i=1}^{n} \cos(c x_i) \right) + a + \epower{},
\end{equation}
where \(\epower{}\) is the Napier constant, approximately equal to \(2.718\), while \( a \), \( b \), and \( c \) are constants, commonly set to the parameters: \(a = 20, \,b = 0.2,\, c = 2\pi\).
The global minimum value of the Ackley function is \(0\) attained at \(x=0\). The values of \(k\) and \(h\) are picked as \(30\) and \(0.2\), respectively. Table \ref{tab:Ackley-2} summarizes the results of \(10\) simulations with various parameters.
%
\begin{table}[htpb]
\centering
\begin{adjustbox}{max width=\textwidth}
\begin{tblr}{lcccccc}
\hline[2pt]
	\SetRow{azure9}
Iterations &  Time taken (seconds) & Log regret & Values attained & Optimal value & Domain \\
\hline[1pt]
200 & 110 & \(\mathbf{0.22997}\) & 1.25856 & 0 & $[-10, 10]^{10^3}$ \\
1000 & 538 & \(\mathbf{-0.490062}\) & 0.61258 & 0 & $[-10, 10]^{10^3}$ \\
\hline[2pt]
\end{tblr}
\end{adjustbox}
\caption{Results for the Ackley function in dimension $10^3$.}
\label{tab:Ackley-2}
\end{table}
Figure \ref{fig:ackband} summarizes the best values obtained with their standard error for \(10\) simulations each consisting of 1000 iterations. Figure \ref{fig:Ack-fig} showcases the values obtained for one such simulation as a function of iterations and the effectiveness of \(\algoname\) in exhibiting fast convergence to the global minimum.
\subsection{Levy function}
\label{subsec:levy function}
\begin{figure}[!htpb]
    \centering
     \includegraphics[width=0.6\linewidth]{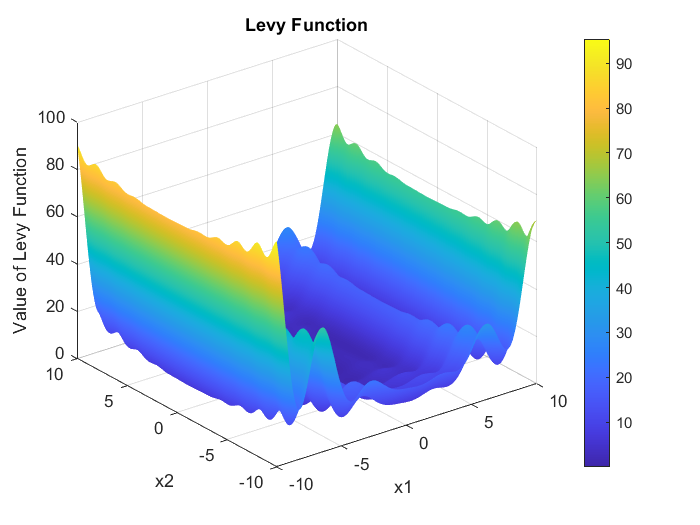} 
     \caption{Plot of the Levy function defined on $[-10,10]^2$.}
    \label{fig:levplot}
\end{figure}

\begin{figure}[htpb]
    \centering
     \includegraphics[width=\linewidth]{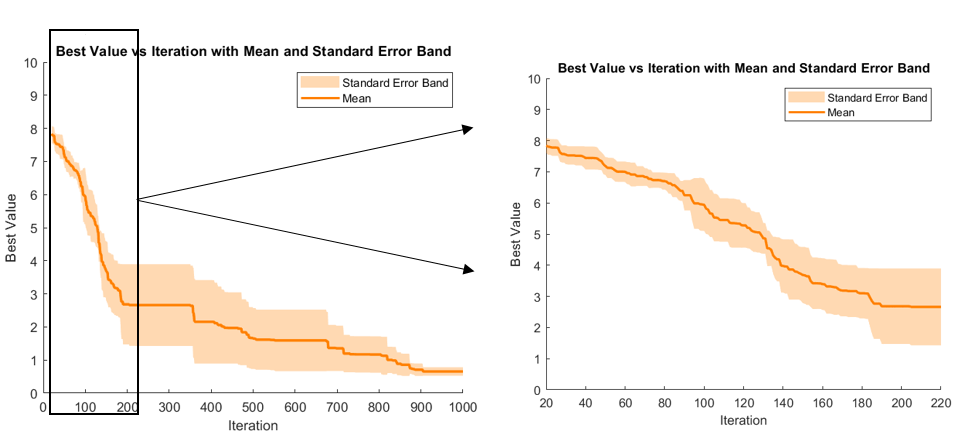} 
     \caption{Plots showing the average log-scaled regret values against iterations for the Levy function in dimension $10^3$. The curve shows the mean value with standard error represented as a band measured across 10 independent simulations. The complete numerical experiment run for 1000 simulations is on the left and a zoomed-in version showing the first 200 iterations out of 1000 is displayed on the right. Here too, the burn-in time is 20 iterations, which is excluded from the figures.}
    \label{fig:Levy-band}
\end{figure}

\begin{figure}[htpb]
    \centering
    \begin{minipage}{0.48\textwidth}
        \centering
        \includegraphics[width=\linewidth]{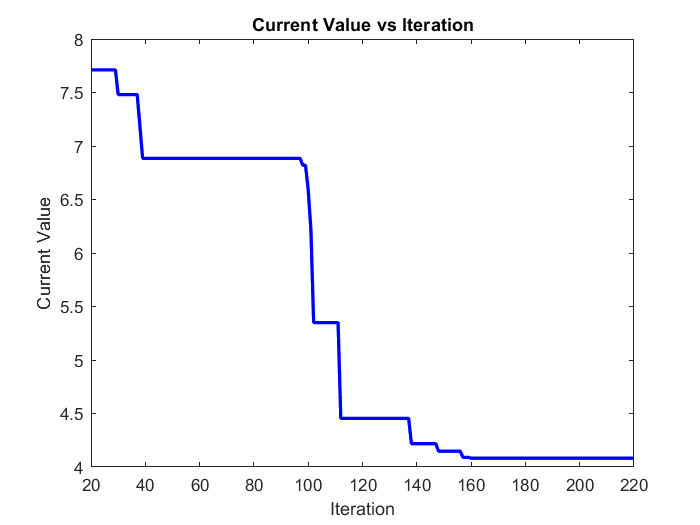} 
        \subcaption{Plot depicting the log regret of the values obtained against iterations.}
        \label{fig:iimage1}
    \end{minipage}
     \hspace{0.02\textwidth}
    \begin{minipage}{0.48\textwidth}
        \centering
        \includegraphics[width=\linewidth]{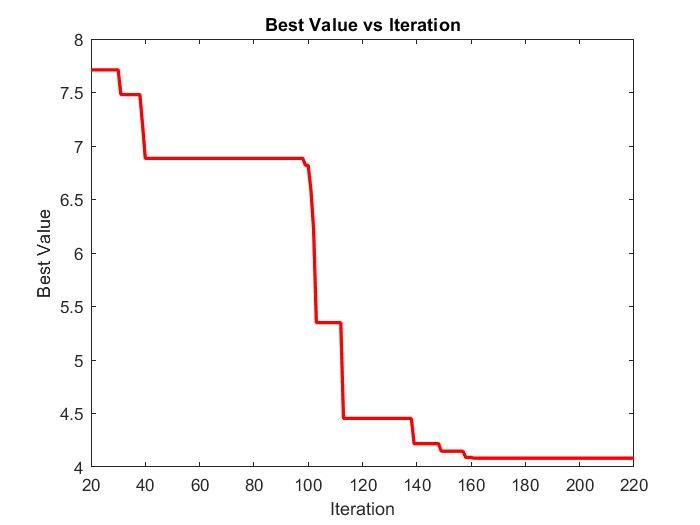} 
        \subcaption{Plot depicting the log regret of the best values obtained thus far against iterations.}
        \label{fig:iimage2}
    \end{minipage}
    \vspace{0.5cm}
    
    \begin{minipage}{0.48\textwidth}
        \centering
        \includegraphics[width=\linewidth]{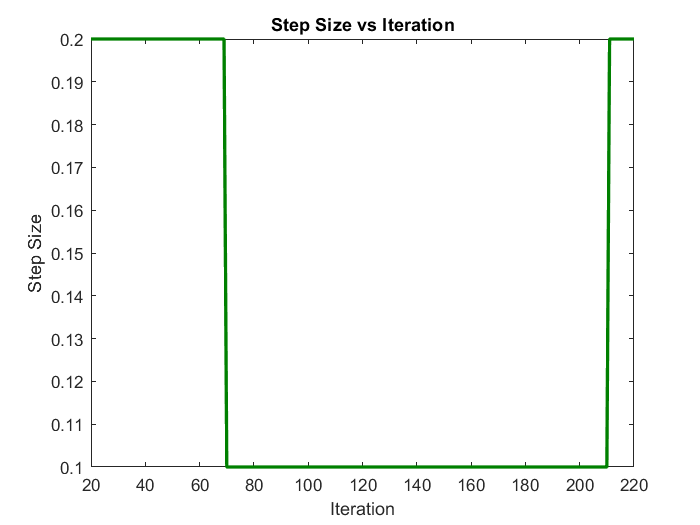} 
        \subcaption{Plot plotting the varying step size trend against iterations.}
        \label{fig:iimage3}
    \end{minipage}
     \hspace{0.02\textwidth}
    \begin{minipage}{0.48\textwidth}
        \centering
        \includegraphics[width=\linewidth]{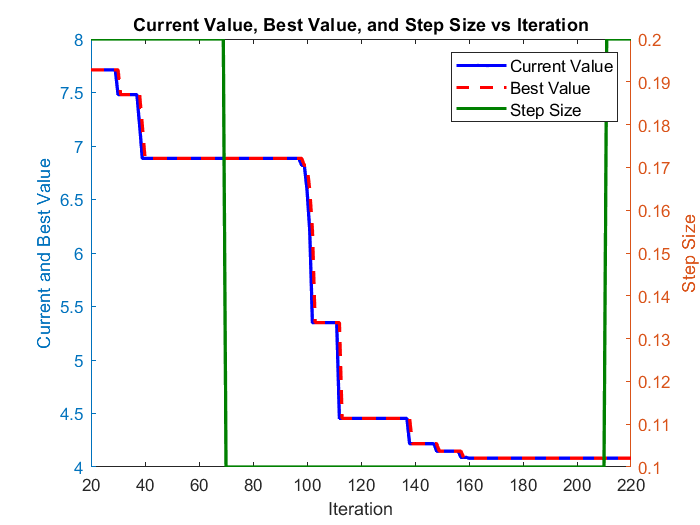} 
        \subcaption{A combined plot showcasing the change in values with the varying step size.}
        \label{fig:iimage4}
    \end{minipage}
    \caption{For the Levy function the simulation was run with a burn-in time of 20 iterations. Here, we plot the results obtained for the next 200 iterations. Graph (d) displays the changes in values of log regret and step size. A finer mesh can be seen to aid the progression towards global minima.}
    \label{fig:Levy-fig}
\end{figure}

The Levy function is another well-known benchmark function for evaluating the performance of optimization algorithms. It is characterized by its numerous local minima and complex, multi-modal landscape. For each \( x_i \) define the auxiliary variable \(w_i\) by
\(
    w_i = 1 + \frac{x_i - 1}{4}.
\)
For an \( n \)-dimensional input vector \( x = (x_1, x_2, \ldots, x_n) \), the Levy function is defined by:
\begin{align}
    \Rbb^n \ni x  \mapsto f(x) = \sin^2(\pi w_1) +& \sum_{i=1}^{n-1} (w_i - 1)^2 \left[ 1  + 10 \sin^2(\pi w_i + 1) \right] \nn \\& + (w_n - 1)^2 \left[ 1 + \sin^2(2\pi w_n) \right]
\end{align}
%
%
Figure \ref{fig:levplot} displays a 2-dimensional plot for the Levy function on $[-10, 10]^2$. 
The global minimum value of the Levy function is \(0\) attained at \(x=1\). \(k\) and \(h\) were chosen to be \(40\) and \(0.2\), respectively. Table \ref{tab:Levy-2} summarizes the results of \(10\) simulations with various parameters. 
%
\begin{table}[!htpb]
\centering
\begin{adjustbox}{max width=\textwidth}
\begin{tblr}{lcccccc}
\hline[2pt]
	\SetRow{azure9}
Iterations & Time taken (seconds) & Log regret & Values attained & Optimal value & Domain \\
\hline[1pt]
200 & 130 & \(\mathbf{2.66177}\) & 14.32161 & 0 & $[-10, 10]^{10^3}$ \\
1000 & 750 & \(\mathbf{0.65386}\) & 1.92294 & 0 & $[-10, 10]^{10^3}$ \\
\hline[2pt]
\end{tblr}
\end{adjustbox}
\caption{Results for the Levy Function in dimension $10^3$.}
\label{tab:Levy-2}
\end{table}
Figure \ref{fig:Levy-band} summarizes the best values obtained with their standard error for \(10\) simulations each consisting of 1000 iterations. Figure \ref{fig:Levy-fig} showcases the values obtained for one such simulation as a function of iterations and the effectiveness of \(\algoname\) in exhibiting fast convergence to the global minimum.
\subsection{Process synthesis MINLP}
The following problem is an example of constrained global optimization lifted from \cite{Gar2023}. 
%
\begin{equation}
    \label{e:arfOCP}
\begin{aligned}
& \minimize && (y_1 - 1)^2 + (y_2 - 2)^2 + (y_3 - 1)^2 - \log(y_4 + 1) \\
&  && \quad + (x_1 - 1)^2 + (x_2 - 2)^2 + (x_3 - 3)^2 \\
& \text{subject to} && \begin{cases}
    y_1 + y_2 + y_3 + x_1 + x_2 + x_3 \leq 5, \\
    y_2^2 + x_1^2 + x_2^2 + x_3^2 \leq 5.5, \\
    y_1 + x_1 \leq 1.2, \,
    y_2 + x_2 \leq 1.8, \\
    y_3 + x_3 \leq 2.5, \,
    y_4 + x_1 \leq 1.2, \\
    y_2^2 + x_2^2 \leq 1.64, \,
    y_3^2 + x_3^2 \leq 4.25, \,
    y_2^2 + x_3^2 \leq 4.64, \\
    (x_1 \, x_2 \,x_3)^{\top} \in \lcrc{0}{1.2} \times \lcrc{0}{1.8}  \times \lcrc{0}{2.5},\\ 
    y_i \in \{0,1\} \text{ for each } i = 1, 2, 3, 4.
\end{cases} 
\end{aligned}
\end{equation}
 Table \ref{tab:minlp} compiles the results obtained.
\begin{table}[htpb]
    \centering
\begin{adjustbox}{max width=\textwidth}
    \begin{tblr}{lcccc}
        \hline[2pt]
		\SetRow{azure9}
        Dimensions & \(k\) & \(h\) & Optimal value obtained & Optimal Value\\
        \hline[1pt]
        7 & 15 & 0.2 & \(\mathbf{4.81}\) & 4.57\\
        \hline[2pt]
    \end{tblr}
\end{adjustbox}
    \caption{Table showcasing the results obtained for the MINLP problem against the reported values lifted from \cite{Gar2023}.}
    \label{tab:minlp}
\end{table}
As the algorithm was built to be integrated with MATLAB's \textsf{simulannealbnd} function that cannot handle constraints efficiently, we did not add them. To overcome this, the constraints were made part of a modified objective function and a heavy penalty was imposed on its value if the constraints were violated. Similar ideas were also used in the next example.

\subsection{Quadratically constrained quartic program}
The next problem is lifted from \cite{Gar2023}.
\begin{equation}
\begin{aligned}
    & \min && x_1^4 - 14x_1^2 + 24x_1 - x_2^2 \\
    & \text{subject to} && \begin{cases}
        -x_1 + x_2 - 8 \leq 0, \\
        x_2 - x_1^2 - 2x_1 + 2 \leq 0, \\
        -8 \leq x_1 \leq 10, \\
        -8 \leq x_2 \leq 10.
    \end{cases}
\end{aligned}
\end{equation}
 Table \ref{tab:quadcons} compiles the results obtained.
\begin{table}[htpb]
    \centering
\begin{adjustbox}{max width=\textwidth}
    \begin{tblr}{lcccc}
        \hline[2pt]
		\SetRow{azure9}
        Dimensions & \(k\) & \(h\) & Optimal value obtained & Optimal value\\
        \hline[1pt]
        2 & 20 & 0.2 & \(\mathbf{-118.168}\) & -118.704\\
        \hline[2pt]
    \end{tblr}
\end{adjustbox}
    \caption{Table showcasing the results obtained for the quadratically constrained quartic problem against the reported values reported in \cite{Gar2023}.}
    \label{tab:quadcons}
\end{table}

\subsection{Chebyshev center problem}
%
A Chebyshev center of a closed and bounded subset \( K \subset \mathcal{Z} \) is defined as the center of a ball of the smallest
radius circumscribing \( K \). In other words, a Chebyshev center of a set \( K \) is an optimizer of the minmax optimization problem:
\begin{equation*}
\inf_{x \in \mathbb{X}} \sup_{y \in K} \| x - y \|
\end{equation*}
This minmax problem is important in the context of learning theory and extracts a single best representative of the hypothesis class that satisfies a given interpolation data; see \cite{ref:PP-DC-23} for details. However, as reported in \cite[Chapter \(15\), pp. no. \(362\)]{ref:AA-IT-21}, even in the simplest case when the decision space \(\mathbb{X}\) and the compact \(K\) are subsets of some Euclidean space, it is NP-hard, and the corresponding time complexity of solving such problems grows exponentially with the dimension of the decision space. Here we solve this problem without assuming the convexity of the set \(K\) and demonstrate that our algorithm \(\algoname\) combined with the MSA algorithm effectively solves the problem near-exactly and in reasonable time. Table \ref{tab:cheby} contains the results obtained.

\begin{figure}[htpb]
\centering
    \includegraphics[width=0.6\linewidth]{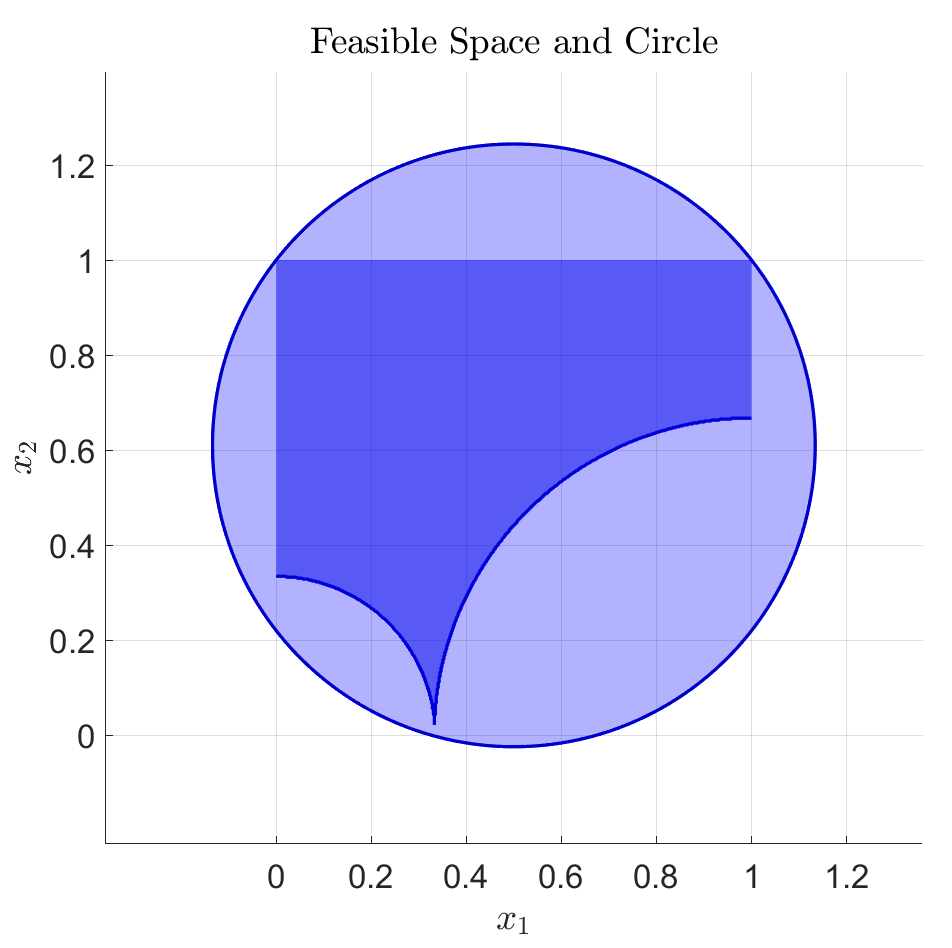}
    \caption{Chebyshev circle and the region \(K\).}
    \label{fig:chebyplot}
\end{figure}

Introducing a slack variable, the Chebyshev center problem can be encoded as the convex semi-infinite program \cite{ref:PP-DC-23}
\begin{equation}
\begin{aligned}
	& \min_{(t, x)} && t \\
    & \text{subject to} && \begin{cases}
        \| x - y \| \leq t \quad \text{for all } y \in K, \\
        (t, x) \in [0, +\infty[ \times \mathbb{Z}.
    \end{cases}
\end{aligned}
\end{equation}
Pick the non-convex set \(K\) from \cite[\S4.5]{ref:PP-DC-23} for our illustration:
\begin{equation*}
K \Let \left\{ (x_1, x_2) \in [0, 1]^2 \, \bigg| \, \sqrt{x_1^2 + x_2^2} \geq \frac{1}{3} \text{ and } \sqrt{(x_1 - 1)^2 + x_2^2} \geq \frac{2}{3} \right\},
\end{equation*}
which consists of a pinched region. The results of \(\algoname\) deployed to solve the preceding problem via the MSA algorithm are reported in Table \ref{tab:cheby}. Figure \ref{fig:chebyplot} displays the non-convex region defined and the Chebyshev circle obtained.
%
\begin{table}[!htpb]
    \centering
\begin{adjustbox}{max width=\textwidth}
    \begin{tblr}{lccccc}
        \hline[2pt]
		\SetRow{azure9}
        \(k\) & \(h\) & Chebyshev radius & Center & Optimal Chebyshev radius & Optimal center\\
        \hline[1pt]
        15 & 0.1 & 0.6341 & (0.4999, 0.6101) & 0.6334 & (0.5000, 0.6111)\\
        \hline[2pt]
    \end{tblr}
\end{adjustbox}
    \caption{Comparison of the results for the Chebyshev center problem against \cite{ref:PP-DC-23}.}
    \label{tab:cheby}
\end{table}

\subsection{SDP - Truss1}
We work with the following SDP written in the standard form used by the SDPA package \cite{ref:BB-99} and encoded in the SDPA sparse format corresponding to the problem Truss \(1\). Table \ref{tab:truss} compiles the results obtained, and the problem itself is:
\begin{equation}
\begin{aligned}
    & \min_{x \in \Rbb^n} && C^\top x \\
    & \text{subject to} && \begin{cases}
        \sum_{i=1}^{m} (x_i F_i - F_0) = X, \\
        X \succeq 0.
    \end{cases}
\end{aligned}
\end{equation}
where \(C \in \mathbb{R}^n\) is a given vector, \(F_i \in \mathbb{R}^{n \times n}\) are given symmetric matrices for \(i = 0, \ldots, m\), and the decision variable \(X \in \mathbb{R}^{n \times n}\) is assumed to be positive semi-definite matrix. 
%
\begin{table}[htpb]
\centering
\begin{adjustbox}{max width=\textwidth}
\begin{tblr}{ccccccc}
\hline[2pt]
	\SetRow{azure9}
\(m\) & \(n\) & \(k\) & \(h\) & Value attained & Optimal value & Domain (x) \\
\hline[1pt]
6 & 13 & 20 & 0.4 & \(\mathbf{-8.84306}\) & -8.99 & [-4, 4] \\
\hline[2pt]
\end{tblr}
\end{adjustbox}
\caption{Table displaying the results of the SDP problem.}
\label{tab:truss}
\end{table}

\section{Conclusion}\label{sec:conclu_dissc}
We introduced a parallelizable probabilistic zeroth-order global optimization algorithm for solving a fairly general class of problems. The key components driving the algorithm are the Walker slice sampling scheme and the Gibbs sampler. The effectiveness of the algorithm was tested on several benchmark problems and comparisons, both in terms of speed and accuracy, against the state-of-the-art algorithms in the literature were reported.

The algorithm was developed for global optimization algorithms that arise as an intermediate step in solving convex semi-infinite programs \cite{ref:DasAraCheCha-22}. These problems have immense applications in various domains of science and engineering. 
We emphasize that \(\algoname\) works under weak hypotheses (in particular, convexity/differentiability/continuity need not be assumed) of measurability on the objective function, it is quick due to the discretization of the state-space, accurate because of the adaptive refinement scheme, scales well with the number of dimensions, and is completely parallelizable, making it suitable for high-dimensional optimization problems.

\bibliographystyle{plain}
\bibliography{ref.bib}

\appendix

\section{An exposition to Walker slice sampling techniques}\label{subsec:ws}
Let \( \pi(\cdot) \) be the stationary distribution of the Markov chain defined on \(\ol{\stsp}\) with the transition density \( p_A\cprobof[]{\cdot \given \cdot} \). Pick the \emph{annealing temperature} \(\temp>0\). We choose \(\pi\) to be the unnormalized Boltzmann distribution given by the following expression: 
\begin{equation}
    \label{eq:pi}
    \pi(r) = \exp{\bigl( -\obj(r) / \temp \bigr)} \quad \text{for }r \in \ol{\stsp}.
\end{equation}
Fix the \emph{exploration parameter} \( k \in \N \). For states \(r,s \in \ol{\stsp}\) such that \( |r - s| \leq k - 1 \), the transition probability mass function of the Walker slice sampler is given by
\begin{equation}\label{eq:transition probability}
    p_A\cprobof[]{r\given s} = \frac{\pi(r)}{k} \sum_{\ell=\max\{s, r\}}^{\min\{s+k-1, r+k-1\}} \frac{1}{\sum_{j=\max\{1, \ell-k+1\}}^{l} \pi(j)}.
\end{equation}
For a fixed choice of the exploration parameter \( k \), the total number of neighboring states that can be explored is  \( 2k - 1 \). The following result is crucial: 
\begin{proposition}[\cite{ref:SGW-14}]
    \label{prop:WS sampler}
    Consider the Walker slice sampler with the transition probability mass function \(p_A \cprobof[]{\cdot \given \cdot }\) defined in \eqref{eq:transition probability}. Denote the normalized version of \( \pi \) by 
    \begin{equation}\label{eq:unnorm pi}
        \widetilde{\pi}(r) = \frac{1}{Z}\pi(r) \quad \text{for }r \in \ol{\stsp},
    \end{equation} 
    where \(Z>0\) is the normalization parameter. If \(r_0 = \max\{1, s-k+1\}\), then \(\sum_{r=r_0}^{s+k-1} p_A\cprobof[]{r\given s} = 1\). Moreover,  
    \begin{align*}
        p_A\cprobof[]{r\given s} \xrightarrow[k \to +\infty]{} \widetilde{\pi}(r) \quad \text{for all \(r,s, \in \ol{\stsp}\)}.
    \end{align*}
\end{proposition}
Proposition \ref{prop:WS sampler} asserts that the transition probability converges to the stationary distribution \(\widetilde{\pi}\); consequently, one can replace the Metropolis-Hasting sampling scheme in the simulated annealing algorithm with the Walker slice sampler. It is straightforward to see that the transition density \( p_A \cprobof[]{\cdot \given \cdot}\) defined in \eqref{eq:transition probability} also satisfies the detailed balance condition
\begin{equation*}
    p_A \cprobof[]{r \given s} \pi(s) = p_A \cprobof[]{s \given r} \pi(r) \quad \text{for every }r,s \in \ol{\stsp}. 
\end{equation*}
Readers are referred to \cite{ref:SGW-14} for more details on this sampling scheme.
\section{Details of the addendum to Algorithm \ref{alg:swiftnav}: \(\mathsf{ProbCal}\) and \(\mathsf{AdapRef}\)}
\label{appen:details}

\subsection{\(\mathsf{ProbCal}\)}
\label{appen:subsec:probcal}
Algorithm \ref{alg:probcal} implements the Walker slice sampling routine. 
Algorithm \ref{alg:probcal} computes the probability \(p_A\cprobof[\big]{r_i(t) \given s_{-i}(t)}\) for the \(i\)$^{\text{th}}$ dimension and it is invoked in Algorithm \ref{alg:swiftnav}. The following notations are used in Algorithm \ref{alg:probcal}:
\begin{enumerate} [leftmargin=*,label=(\Alph*)]
    


    
    \item \label{it:proObj}The variable \(\text{\texttt{proObj}}\) stores the value of \(s(t)\). Note that the \(i\)$^{\text{th}}$ component of \(\text{\texttt{proObj}}\) is modified keeping all other dimensions unchanged.

    \item Let \(\text{\texttt{lb}}\) and \(\text{\texttt{ub}}\) denote the lower bound and upper bound, respectively, of the feasible region of the optimization problem \eqref{eq:global opt prob} along each dimension. 

    
    
    \item  \(\text{\texttt{pi(\text{\ttt{r}})}}\) evaluates \(\pi\) at \(\text{\texttt{r}}\) and is used to calculate \( p_A \cprobof[]{\cdot \given \cdot}\). The values \(p_A \cprobof[]{\text{\ttt{r}} \given \text{\ttt{s}}}\) are stored in \(\text{\texttt{pA}}\).    
\end{enumerate}
\begin{algorithm}[ht]
\footnotesize
\SetAlgoLined
\SetKwInOut{Input}{Input}
\SetKwInOut{Output}{Output}
\SetKwInOut{ini}{Initialize}
\Input{\(\text{\texttt{r}}\),  \(\text{\texttt{s}}\), \(k\), \(\text{\texttt{currdim}}\), \(\text{\texttt{neigh}}\), \(\text{\texttt{currentState}}\), \(\text{\texttt{lb}}\), \(\text{\texttt{ub}}\), \(\text{\texttt{objectiveFunc}}\), \(\text{\texttt{optimValues}}\), \(\text{\texttt{pi\_table}}\)}
\ini{\(\text{\texttt{l\_min}} \gets \max(\text{\texttt{s}}, \text{\texttt{r}})\), \(\text{\texttt{l\_max}} \gets \min(\text{\texttt{s}} + k - 1, \text{\texttt{r}} + k - 1)\)}
\Output{Probability \(\text{\texttt{pA}}\) and updated table \(\text{\texttt{pi\_table}}\)}

\SetKwFunction{Fpi}{pi}
\SetKwProg{Fn}{Function}{:}{end}
\Fn{\Fpi{r}}{
    \If{\(\text{\texttt{pi\_table(r)}} == -1\)}{
        \(\text{\texttt{proObj}} \gets \text{\texttt{currentState}}\)\hspace{2mm}\tcp{\(\text{\texttt{currentState}}\) denotes \(s(t)\)}
        \(\text{\texttt{proObj(currdim,1)}} \gets \text{\texttt{neigh(r,1)}}\) \hspace{2mm} \tcp{\text{\texttt{currdim}} is the index of the dimension. Vector containing all the \(2k-1\) neighbors of \(s\) is represented by  \(\text{\ttt{neigh}}\)}
        \If{\(\text{\texttt{proObj(currdim,1)}} < \text{\texttt{lb(currdim,1)}} \,||\, \text{\texttt{proObj(currdim,1)}} > \text{\texttt{ub(currdim,1)}}\)}{
            \(\text{\texttt{pie}} \gets 0\)\;
        } \hspace{2mm} \tcp{\(\text{\texttt{lb}}\) and \(\text{\texttt{ub}}\) denote the lower bound and upper bound, respectively, of the feasible region of \eqref{eq:global opt prob}}
        \Else{
            \(\text{\texttt{pie}} \gets \exp{\biggl( \frac{-\text{\texttt{objectiveFunc(proObj)}}}{\text{\texttt{optimValues.temperature(currdim,1)}}}\biggr)}\)\;
        }
        \(\text{\texttt{pi\_table(r)}}\) $\gets$ \(\text{\texttt{pie}}\)\;
        }
    
    \Else{
        \(\text{\texttt{pie}}\) $\gets$ \(\text{\texttt{pi\_table(r)}}\)\;
    }
}
\(\texttt{sum\_l} \gets 0\) \;

\For{\(l = \texttt{l\_min} : \texttt{l\_max}\)}{
    \texttt{sum\_j} $\gets$ 0\;
	
    \For{\(j = \max(1, l - k + 1) : \texttt{l}\)}{
        \texttt{sum\_j} += \Fpi{j}\;
    }
    \texttt{sum\_l} += $\frac{1}{\texttt{sum\_j}}$\;
}

Return \(\text{\texttt{pA}} \gets \frac{\Fpi{\texttt{r}} \cdot \texttt{sum\_l}}{k} \)\;
\caption{\textsf{ProbCal}: Routine to compute \(p_A \cprobof[]{\text{\ttt{r}} \, \mid \, \text{\ttt{s}}}\)}
\label{alg:probcal}
\end{algorithm}

\subsection{\(\mathsf{AdapRef}\)}
\label{appen:subsec_adapref}
 Algorithm \ref{alg:adapref} determines if the mesh-refinement scheme is necessary to enhance the accuracy of the obtained optimal solution by fine-tuning the step size \(h\) provided that certain conditions are satisfied. Below we fix the notations in Algorithm \ref{alg:adapref}:
%
%
%
\begin{enumerate} [leftmargin=*,label=(\Alph*)]
    

    \item \(\texttt{iterationCount}\) stores the number of iterations passed since Algorithm \ref{alg:adapref} received a value smaller than the \(\texttt{bestValue}\). It resets whenever the condition mentioned above is satisfied.

    \item Denote by \(\meshfac>0\) the factor by which \(h\) is reduced. Let \(p>0\) be chosen in an ad-hoc manner. Suppose that 
    \(
        \texttt{reducedStepSizeFlag} \in \aset[]{0,1},
    \) 
    and the boolean value admitted by \(\texttt{reducedStepSizeFlag}\) activates the adaptive-refinement scheme. We have the following two cases:
    \begin{itemize}[label = \(\circ\)]
        \item Define the interval 
        \[
        \mathbb{I}   \Let \lcrc{\texttt{bestValue} - 0.1 \texttt{bestValue}}{\texttt{bestValue} + 0.1\texttt{bestValue}}.
        \]
        \(\texttt{reducedStepSizeFlag}\) is set to \(1\) only when \(\texttt{currentValue} \in \mathbb{I}\) and \(\texttt{iterationCount}>p\), reducing the step size \(h\) by \(\mathtt{\meshfac}\).

        \item \(\texttt{reducedStepSizeFlag}\) is set to \(0\). It refers to the adaptive-refinement scheme being not active with the default value of 0.
    \end{itemize}

    \item \(\texttt{reducedIterationCount}\) functions similar to \(\texttt{iterationCount}\). When $\texttt{reducedStepSizeFlag}$ is set to \(1\), \(\texttt{reducedIterationCount}\)  stores the number of iterations that have passed since the Algorithm \ref{alg:adapref} has received a value smaller than the \texttt{bestValue}.
    
    \item The mesh size can be further reduced if \texttt{currentValue} improves further by checking if \texttt{reducedIterationCount} exceeds a limit \(q>0\) chosen beforehand. The step size \(h\) is then further refined by \(\meshfac\). If no improvement is observed, \(h\) is reset to its initial value. \texttt{reducedIterationCount} is reset in both cases.  
\end{enumerate}
\begin{algorithm}[!htpb]
\footnotesize
\SetAlgoLined
\SetKwInOut{Input}{Input}
\SetKwInOut{Initialize}{Initialize}
\SetKwInOut{Output}{Ouput}

 \Input{\(\texttt{currentValue}\), \( \texttt{bestValue}\), \(\) \texttt{iterationCount}, \(\texttt{reducedStepSizeFlag}\), \(\texttt{reducedIterationCount}\), \(p\), \(q\), \(\meshfac\)}
\If{\(\texttt{currentValue} < \texttt{bestValue}\) }{
    \(\texttt{oldBenchmark} \gets \texttt{bestValue}\) \tcp{The variable \(\texttt{bestValue}\) stores the minimum value evaluated by running Algorithm \ref{alg:swiftnav}, until time \(t\). \texttt{currentValue} is the evaluation of \(\text{\texttt{objectiveFunc}}\)at \(s(t+1)\) in Algorithm \ref{alg:swiftnav}; see Step \ref{line:obj_init} in Algorithm \ref{alg:swiftnav}} 
    \(\texttt{iterationCount} \gets 0\)\;
    \(\texttt{bestValue} \gets \texttt{currentValue}\) \;
    \If{\(\texttt{reducedStepSizeFlag} == 1\)}{
        \(\texttt{reducedIterationCount} \gets 0\)\;
    }
}
\Else{
    \If{\(\texttt{reducedStepSizeFlag} == 1\)}{
        \(\texttt{reducedIterationCount}  \gets \texttt{reducedIterationCount} + 1\)\;
    }
    \Else{
        \(\texttt{iterationCount} \gets \texttt{iterationCount} + 1\)\;
    }
}
\If{\(\texttt{iterationCount} >= p \; \&\&\; \texttt{currentValue} \in \mathbb{I}\)}{
    Reduce step size \(h \gets \frac{h}{\mathtt{\meshfac}}\)\;
    \(\texttt{reducedStepSizeFlag} \gets 1\)\;
    \(\texttt{reducedIterationCount} \gets 0\) \;
}
\If{\(\texttt{reducedStepSizeFlag} == 1 \; \&\& \; \texttt{reducedIterationCount} >=q\)}{
    \If{\(\texttt{oldBenchmark} \in \mathbb{I}\)}{
        Restore step size \(h\)\;
        \(\texttt{reducedStepSizeFlag} \gets 0\) \;
        \(\texttt{reducedIterationCount} \gets 0\)\;
    }
    \Else{
        Reduce step size \(h \gets \frac{h}{\mathtt{\meshfac}}\)\;
        \(\texttt{reducedIterationCount} \gets 0\)\;
    }
}
\caption{\ttt{AdapRef}: Routine to implement the adaptive-refinement scheme}
\label{alg:adapref}
\end{algorithm}
\section{Background}
\label{appen:background}

\subsection*{Simulated Annealing}
Here we provide a brief description of the simulated annealing algorithm; interested readers
are referred to \cite{Kirkpatrick1984}, \cite{Cerny1985} for a detailed exposition. The simulated
annealing algorithm is an adaptation of the Metropolis-Hastings algorithm \cite{Robert1999}, \cite{Chib1995} to obtain the global optimum of a given function.
Suppose that we intend to obtain the minimum of a function \( \obj : V \rightarrow [0, +\infty[ \), where
\( V \subset \mathbb{R}^k \). Let \( V \ni \zeta \rightarrow \widetilde{\pi}(\zeta) := \frac{1}{Z} e^{\obj(\zeta)/\temp} \in [0, +\infty[ \) be a density on \( \mathbb{R}^k \) where \( Z \) is a
normalization constant. Here, the parameter \( \temp > 0 \) is the control parameter. It is analogous to
temperature in annealing (metallurgy) and as the physical temperature behaves in annealing,
the parameter \( \temp \) is reduced with the progress of the simulated annealing algorithm. To reduce
\( \temp \), a cooling sequence denoted by \( (\temp_N)_{N \geq 1} \) is used, where \( N \) denotes the current iteration
number of the algorithm. Note that since \( \widetilde{\pi}(\zeta) \propto e^{-\obj(\zeta)/\temp} \), the density \( \widetilde{\pi} \) has its maximum
value at the global minimum of \( \obj \). Thus the samples drawn from \( \widetilde{\pi} \) are concentrated near
the global minimum value of \( \obj \). The simulated annealing algorithm utilizes this idea, along
with a modification of the Metropolis-Hastings algorithm, to obtain the global minimum of
\( \obj \).

\subsection*{Gibbs sampler}
Gibbs Sampling \cite{ref:Hagg-02} is an MCMC algorithm for generating samples from a joint probability distribution when direct sampling is challenging. It leverages conditional distributions to update each variable iteratively. Below we enumerate the key steps:\\

\textbf{Algorithm description:}
Given a random vector \( {X} = (X_1, X_2, \ldots, X_n) \) with joint distribution \( P({X}) \), where \(n\) represents the number of dimensions and \( X_i \) comprises of the current values; the Gibbs sampler is a Markov chain which proceeds as follows:
\begin{enumerate}
    \item \textbf{Initialization}: Start with a random initial state \( X^{(0)} \).
    \item \textbf{Iteration}: For iteration \( t = 1, 2, \ldots \), update each component \( X_i \) according to their conditional probability distribution given that all the other dimension values remain unchanged, i.e.,
     \begin{align*}
         \begin{cases}
             X_1^{(t+1)} \sim \PP\cprobof[\big]{X_1 \given X_2^{(t)}, X_3^{(t)}, \ldots, X_n^{(t)}} \\
             X_2^{(t+1)} \sim \PP\cprobof[\big]{X_2 \given X_1^{(t+1)}, X_3^{(t)}, \ldots, X_n^{(t)}}\\
            \hspace{2cm}\vdots\\
             X_n^{(t+1)} \sim \PP \cprobof[\big]{X_n \given X_1^{(t+1)}, X_2^{(t+1)}, \ldots, X_{n-1}^{(t+1)}}.
         \end{cases}
     \end{align*}

    \item After sufficient iterations, the samples \(X^{(t)} \) approximate the target distribution \( \PP\probof[\big]{X}\).
\end{enumerate}

\end{document}